\documentclass[12pt]{amsart}
\usepackage{amscd,amssymb}
\usepackage[graph,frame,poly,arc]{xy}  
\usepackage[plainpages,backref,urlcolor=blue]{hyperref}

\theoremstyle{plain}
\newtheorem{thm}[subsection]{Theorem}
\newtheorem{lem}[subsection]{Lemma}
\newtheorem{prop}[subsection]{Proposition}
\newtheorem{cor}[subsection]{Corollary}

\theoremstyle{definition}
\newtheorem{rk}[subsection]{Remark}
\newtheorem{definition}[subsection]{Definition}
\newtheorem{ex}[subsection]{Example}

\newtheorem{question}[subsection]{Question}

\numberwithin{equation}{section}
\newcommand{\OO}{{\mathcal O}}

\newcommand{\A}{{\mathcal A}}

\newcommand{\J}{{\mathcal J}}
\newcommand{\D}{{\mathcal D}}

\newcommand{\Z}{\mathbb{Z}}
\newcommand{\Q}{\mathbb{Q}}

\newcommand{\C}{\mathbb{C}}
\newcommand{\PP}{\mathbb{P}}

\newcommand{\N}{\mathbb{N}}

\DeclareMathOperator{\rank}{rank}

\DeclareMathOperator{\pd}{pd}

\DeclareMathOperator{\defect}{def}

\DeclareMathOperator{\depth}{depth }
\DeclareMathOperator{\reg}{reg}

\begin{document}

\title [Free and nearly free surfaces in $\PP^3$]
{Free and nearly free surfaces in $\PP^3$}

\author[Alexandru Dimca]{Alexandru Dimca$^1$}
\address{Universit\'e C\^ ote d'Azur, CNRS, LJAD, France }
\email{dimca@unice.fr}

\author[Gabriel Sticlaru]{Gabriel Sticlaru}
\address{Faculty of Mathematics and Informatics,
Ovidius University,
Bd. Mamaia 124, 900527 Constanta,
Romania}
\email{gabrielsticlaru@yahoo.com }
\thanks{$^1$ Partially supported by Institut Universitaire de France.}

\subjclass[2010]{Primary 14J70; Secondary  13D02, 13P20, 14C20, 32S22}

\keywords{Jacobian ideal, Milnor algebra, free divisor, nearly free divisor, Saito's criterion}

\begin{abstract} We define the nearly free surfaces in $\PP^3$ and show that the Hilbert polynomial of the Milnor algebra of a free or nearly free surface in $\PP^3$ can be expressed in terms of the exponents. An analog of Saito's criterion of freeness in the case of nearly free divisors is proven and examples of irreducible free and nearly free surfaces 
are given.
\end{abstract}
 
\maketitle


\section{Introduction} 

Let $S=\oplus_k S_k=\C[x_0,...,x_n]$ be the graded polynomial ring in $n+1$ indeterminates with complex coefficients, where $S_k$ denotes the vector space of degree $k$ homogeneous polynomials. Consider  for a degree $d$ polynomial $f \in S_d$, the corresponding Jacobian ideal $J_f$ generated by the partial derivatives $f_j$ of $f$ with respect to $x_j$ for $j=0,...,n$ and  the graded Milnor algebra $M(f)=\oplus_k M(f)_k=S/J_f$.

The Hilbert function $H(M(f))$ and the Hilbert polynomial $P(M(f))$ of the graded $S$-module $M(f)$ encode information on the 
 projective hypersurface $V=V(f):f=0$  in $\PP^n$ and the associated singular subscheme 
 $\Sigma(f)$ defined by the Jacobian ideal. As an example, when $V$ has only isolated singularities, then the Hilbert polynomial $P(M(f))$ is a constant, equal to the total Tjurina number of $V$, see \cite{CD}, which is also the degree of the singular subscheme $\Sigma(f)$.
On the other hand, the study of the Hilbert function $H(M(f))$ tells much about the syzygies among the partial derivatives $f_i$ and about the geometry of $V$ and of its complement, 
see \cite{DBull, DS1, DS14, DStNach, DStEdin, Se}.

A related subject is the study of free divisors, started by Kyoji Saito \cite{KS} and attracting much interest ever since. In a recent paper \cite{DStNF}, we have  introduced the class of nearly free curves in $\PP^2$, in an attempt to clarify the relation between the rational cuspidal curves and the free curves in the plane, see Remark \ref{rkB+} below for more details. 
In the case of surfaces considered in this note, we  follow essentially the same definition as in the curve case, in the hope that this  leads to an interesting class of surfaces. An unexpected byproduct is the relation between free and nearly free surfaces and homaloidal polynomials, see \cite{DStMM}.

Since the free and the nearly free surfaces in $\PP^3$  have a 1-dimensional singular locus, in the second section we collect some basic facts on the Hilbert function and the Hilbert polynomial of a projective hypersurface $V=V(f):f=0$  in $\PP^n$  for
$n \geq 3$ whose singular locus $\Sigma$ satisfies $\dim \Sigma=1$. To illustrate the results stated there, we consider in the third section two simple cases: the case when $V$ is a surface obtained as the union $D\cup D'$ of two smooth surfaces $D$ and $D'$ in $\PP^3$ meeting transversally, and the case when $V$ is a cone (in two natural ways) over a hypersurface $W$ in $\PP^{n-1}$ having only isolated singularities.
The latter case can be used to show by examples that the Hilbert polynomial $P(M(f))$ may depend on the position of the singularities, e.g. when the singular locus $\Sigma$ consists of three concurrent lines, then the polynomial $P(M(f))$ may depend on whether or not these lines are coplanar, see 
Examples \ref{ex2} and  \ref{ex2.5}. It also allows the construction of free (resp. nearly free) surfaces in $\PP^3$ as cones over free (resp. nearly free) curves in $\PP^2$, see Corollaries \ref{corFC} and  \ref{corFC2}, and Definition \ref{NFdef}.

The last two sections contain the main results of this note. Theorem \ref{thmFREE} (resp. Theorem \ref{thmNFREE}) express the Hilbert polynomial $P(M(f))$ and other invariants of the free (resp. nearly free) surface $D:f=0$ in terms of the exponents $d_1\leq d_2 \leq d_3$ of $D$.
In fact, for a free surface, the exponents and the Hilbert polynomial $P(M(f))$ determine each other, see Corollary \ref{corFREE}, (i).

In Theorem \ref{thmSaito} we prove an analog of Saito's criterion of freeness in the case of nearly free surfaces, which expresses the (unique) second order syzygy in terms of determinants constructed using the first order syzygies. A similar result holds for the nearly free curves in $\PP^2$ with the same proof as for surfaces, but was not stated in \cite{DStNF}.

Exactly as in the case of free curves in $\PP^2$ discussed in \cite{ST}, \cite{DStFD}, the irreducible free (resp. nearly free) surfaces are not easy to find. We give both isolated cases and countable families of examples of such surfaces in Examples \ref{exFREE2} and \ref{exNF}. 
More involved examples, related to the discriminants of binary forms, are given in Propositions 
\ref{propD3}, \ref{propD4} and \ref{propD4'}.
All our examples are rational surfaces (either one of the variables $x,y,z,w$ occurs only with exponent 1 in the defining equations, or as in Propositions  \ref{propD3} and \ref{propD4'} this follows from the description of the discriminants).

Theorem \ref{thmNF} and Example \ref{exNF4}  discuss  the point whether for a nearly free surface, the first local cohomology group $H^1_Q(M(f))$ of the Milnor algebra $M(f)$ with respect to the maximal homogeneous ideal $Q$ in $S$ is a finite dimensional $\C$-vector space. By contrast, note that $H^0_Q(M(f))$ is finite dimensional for any hypersurface $V:f=0$ in $\PP^n$, see Remark \ref{rk0dim}. Examples provided by Aldo Conca are listed in Example \ref{exNF4} (ii) and allow one to construct  rank 3 vector bundles on $\PP^3$ which are not direct sum of line bundles, see Remark \ref{locfree}.

Note that a rational cuspidal curve  can be characterized either as a rational curve which is simply-connected, or as an irreducible curve which is homeomorphic to $\PP^1$. By analogy to the case of rational cuspidal curves, we can ask the following. 

\begin{question}
\label{qInt2}
Which geometric conditions on an irreducible  surface $D$ in $\PP^3$  imply that $D$ is either free or nearly free ?
\end{question}

The computations of various minimal resolutions given in this paper were made using two computer algebra systems, namely CoCoA \cite{Co} and Singular \cite{Sing}.
The corresponding codes are available on request, some of them being available in \cite{St}.

\bigskip

We would like to thank the referees for their careful reading of our manuscript and for their suggestions which greatly improved the presentation of our results. Many thanks also to Aldo Conca for useful discussions and the examples recorded in Example \ref{exNF4} below.

\section{Hilbert functions and Hilbert polynomials} 

The {\it Hilbert function} $H(M(f)): \N \to \N$ of the graded $S$-module $M(f)$ is defined by
\begin{equation}
\label{Hfunc}
 H(M(f))(k)= \dim M(f)_k,
\end{equation}
and it is often encoded in the {\it Hilbert-Poincar\'e series} of $M(f)$
\begin{equation}
\label{Hseries}
 HP(M(f);t)= \sum_k\dim M(f)_kt^k.
\end{equation}

It is known that there is a unique polynomial $P(M(f))(t) \in \Q[t]$, called the {\it Hilbert polynomial} of $M(f)$, and an integer $k_0\in \N$ such that
\begin{equation}
\label{Hpoly}
 H(M(f))(k)= P(M(f))(k)
\end{equation}
for all $k \geq k_0$. In analogy to the case of projective hypersurfaces with isolated singularities considered in \cite{DStEdin}, we introduce the {\it stability threshold} 
$$st(V)=st(f)=\min \{q~~:~~  H(M(f))(k)= P(M(f))(k) \text{ for all } k \geq q\}.$$

Let $\Sigma= \Sigma_1 \cup \Sigma_0$, where $\Sigma_1$ denotes the components of the sigular locus $\Sigma$ of dimension 1 as well as  their embedded 0-dimensional components, and $\Sigma_0$ denotes the union of the isolated points in $\Sigma$, with their multiple structure.
Then the points in $\Sigma_0$ corresponds to the isolated singularities of $V$ and their contribution to the Hilbert polynomial  $P(M(f))$ is again a constant $\tau_0(V)$, the total Tjurina number of $V$, which is just the sum of the individual Tjurina numbers of the isolated singularities of $V$.
This fact comes from the well known relation
\begin{equation}
\label{Hpoly+Euler}
 P(M(f))(k)= \chi(\Sigma, \OO_{\Sigma}(k)),
\end{equation}
for any $k \in \Z$, see \cite{Eis}, p. 197. Indeed, $\OO_{\Sigma}=\tilde M(f)$ is the coherent sheaf on $\PP^n$ associated to the graded $S$-module $M(f)$. Moreover, we clearly have
\begin{equation}
\label{Hpoly+Euler2}
 \chi(\Sigma, \OO_{\Sigma}(k))= \chi(\Sigma_1, \OO_{\Sigma_1}(k))+\chi(\Sigma_0, \OO_{\Sigma_0}(k)),
\end{equation}
and $\chi(\Sigma_0, \OO_{\Sigma_0}(k))=\tau_0(V)$ for any $k$.

The general theory of Hilbert polynomials says that the degree of  $P(M(f))$ is given by the dimension of the support of  $\OO_{\Sigma}=\tilde M(f)$. Hence the assumption $\dim \Sigma=1$ implies that
$$P(M(f))(k)=ak+b$$
where 
\begin{equation}
\label{abformula1}
a=\deg(\Sigma_1) \text { and } b= \chi(\Sigma_1, \OO_{\Sigma_1})+\tau_0(V).
\end{equation}
However, the calculation of $a$ and $b$ in general can be very difficult, due to the multiple structure and/or the singularities of the subscheme $\Sigma_1$.

To compute the Hilbert function $H(M(f))$, one can use the general relation between Hilbert functions and Hilbert polynomials, namely

\begin{equation}
\label{Hfunction}
H(M(f))(k)=P(M(f))(k)+\sum_{i=0}^{2}(-1)^i\dim H^i_Q(M(f))_k,
\end{equation}
where $Q=(x_0,...,x_n)$ is the maximal ideal in the graded ring $S$, see \cite{Eis}, Cor. A1.15. To deal with the local cohomology groups $H^i_Q(M(f))$, we recall the following, see \cite{Eis}, Cor. A1.12.

\begin{prop}
\label{prop1}

\noindent (i) There is an exact sequence of graded $S$-modules
$$0 \to H^0_Q(M(f)) \to M(f) \to \sum_kH^0(\Sigma, \OO_{\Sigma}(k))\to H^1_Q(M(f)) \to 0.$$
\noindent (ii) For every $i \geq 2$,
$$H^i_Q(M(f)) = \sum_kH^{i-1}(\Sigma, \OO_{\Sigma}(k)).$$
In particular $H^i_Q(M(f)) =0$ for $i>2$ since $\dim \Sigma=1$.
\end{prop}
Let $I_f$ denote the saturation of the Jacobian ideal $J_f$ with respect to the ideal $Q$. Then it is clear that 
\begin{equation}
\label{sat}
 H^0_Q(M(f))=I_f/J_f.
\end{equation}
We denote this quotient $I_f/J_f=H^0_Q(M(f))$ that occurs quite frequently by $N(f)$, and the first cohomology group $H^1_Q(M(f))$ by $P(f)$. The exact sequence in (i) above implies
\begin{equation}
\label{ev}
 0 \to S_k/I_{f,k} \to H^0(\Sigma, \OO_{\Sigma}(k)) \to P(f)_k \to 0.
\end{equation}
The arrow $ev:S_k/I_{f,k} \to H^0(\Sigma, \OO_{\Sigma}(k))$ can be thought of as an evaluation map, and by analogy to the case $\dim \Sigma=0$, we define the defect of $\Sigma$ with respect to degree $k$ polynomials to be
$$\defect_k \Sigma= \dim P(f)_k.$$
With this notation, the formula \eqref{Hfunction} becomes
\begin{equation}
\label{Hfunction2}
H(M(f))(k)=P(M(f))(k)+  \dim N(f)_k -\defect_k \Sigma +\dim H^1(\Sigma, \OO_{\Sigma}(k)).
\end{equation}

\begin{rk}
\label{rk0dim} (i) 
In the case when $\Sigma$ is 0-dimensional, one has the similar formula
$$H(M(f))(k)=\tau(V)+\dim N(f)_k -\defect_k \Sigma, $$
see for instance \cite{DBull} or the formula (3.4) in \cite{DStNF}. This is a special case of the formula \ref{Hfunction2}, since clearly $H^1(\Sigma, \OO_{\Sigma}(k))=0$ when $\Sigma$ is 0-dimensional.

\noindent (ii) Note that $H^0_Q(M(f))$ is a finite dimensional $\C$-vector space, as it follows from the general discussion in  \cite{Eis}, pp. 187-188. On the other hand, $H^1_Q(M(f))$ is often 
an infinite dimensional $\C$-vector space, see for instance \cite{Eis}, Example A1.10 and Example \ref{exNF4} below.
\end{rk}

\section{Some simple cases: transversal intersections and cones} 

\subsection{Transversal intersection of two surfaces in $\PP^3$} \label{normalcross}

In this section we consider two smooth surfaces: $D:g=0$ of degree $e$ and $D':g'=0$ of degree $e'$, meeting transversally along the curve $C:g=g'=0$. Since $C$ is a smooth complete intersection, its genus is given by the formula
\begin{equation}
\label{genus}
g(C)=1+(e+e'-4)ee'/2,
\end{equation}
see \cite{D1}, p. 152. Moreover, the degree of $C$ is $ee'$, the number of intersections of $C$ with a generic plane.

Consider now the surface $V=D \cup D': f=gg'=0$. Then the singular scheme $\Sigma$ is reduced, see Proposition \ref{prop2} below, and coincides with the smooth curve $C$. Using \eqref{abformula1}, we get
$$a= \deg(C)=ee' \text {  and  } b=\chi (C,\OO_C)=1-g=-(e+e'-4)ee'/2.$$
Numerical examples suggest the following.
\begin{question}
\label{q0}
Prove that the stability threshhold $st(V)$ in this case is given by the formula
$$st(V)=3(e+e')+ |e-e'|-7.$$
\end{question}

We have also the following result, which actually holds for the transverse intersection of two smooth hypersurfaces in $\PP^n$ for arbitrary $n$ with the same proof as that given below.

\begin{prop}
\label{prop2}
The ideal $I_f$ coincide with the ideal $I(g,g')$ spanned by the two polynomials $g$ and $g'$ in $S$.

\end{prop}

\proof

Since $f=gg'$, it follows that the partial derivatives $f_j=g_jg'+gg'_j$ are in the ideal $I(g,g')$.
Hence $I_f=sat(J_f) \subset sat(I(g,g'))=I(g,g')$, the last equality following from the fact that the ideal $I(g,g')$ is a complete intersection, see \cite{DBull}, Proposition 1.

To prove the converse inclusion, note that any element $h \in J_f$ can be written as
$$h=\sum_ja_jf_j=(\sum_ja_jg_j)g'+(\sum_ja_jg'_j)g.$$
If we take $a_0=g'_1$ and $a_1=-g'_0$ and set $a_j=0$ for $j>1$, then we get that
$$(g_0g'_1-g_1g'_0)g' \in J_f.$$
If $m_2(g,g')$ denotes the ideal in $S$ spanned by all the $2 \times 2$ minors of the matrix formed by the partial derivatives of $g$ and $g'$, then the above shows that $m_2(g,g')g' \subset J_f$.

Since clearly $f=gg' \in J_f$, it is enough to show that the ideal $m_2(g,g')+ (g)$ contains a power $Q^s$ of the maximal ideal $Q$. This is equivalent to showing that the zero set $Z=Z(m_2(g,g')+ (g))$ consists only of the origin. Let $z \in Z$ be different from $0$. Then two cases can occur. If $g'(z)=0$ it follows that $z$ corresponds to a point in the intersection $C$ of the two surfaces, and this is a contradiction, since then at least one minor should not vanish at $z$.
If $g'(z) \ne 0$, then we get that $g(z) \ne 0$. Indeed, the two differentials $dg(z)$ and $dg'(z)$ are non-zero ($D$ and $D'$ are smooth), they are proportional and using the Euler formula for $g(z)$ and $g'(z)$ we get the claim. But this is again a contradiction, since we have supposed $g(z)=0$.

\endproof

\begin{rk}
\label{rkselfdual}
Even in this very simple situation, the graded $S$-module $N(f)$ is no longer self-dual as was the case when $\Sigma$ is 0-dimensional, see \cite{DS1}. For instance, the Hilbert-Poincar\'e series  for $N(f)$ in the case $e= e'= 3$ is given by
$$HP(N(f);t)=2t^3+8t^4+16t^5+23t^6+26t^7+22t^8+12t^9+3t^{10}.$$
\end{rk}

\begin{ex}
\label{ex1}
Consider the case $e=2$ and $e'=3$. Then if the two surfaces $D$ and $D'$ are smooth and meet transversely, e.g. $g= x^2+3y^2+5z^2+7w^2$ and $g'=x^3+y^3+z^3+w^3$, then the Hilbert polynomial has the form
$$P(M(f))(k)=ak+b=6k-3,$$
using the above formulas for $a$ and $b$.
Now consider the more general case when the intersection is transversal, i.e. $C:g=g'=0$ is a smooth complete intersection, but $D$ and $D'$ are allowed to have isolated singularities outside $C$. Then the new formulas for $a$ and $b$ are
$$ a= \deg C=ee'  \text{   and   }  b=\chi (C,\OO_C)+ \tau_0(V)=-(e+e'-4)ee'/2+ \tau_0(V),$$
as in \ref{abformula1}. To have an example, let $g= x^2+3y^2+5z^2+7w^2$ and $g'=xyz+w(xy+yz+xz)$. Then $D$ is a smooth conic, $D'$ is a cubic surface with four $A_1$ singularities and hence $\tau_0(V)=4$. It follows that in this case
$$P(M(f))(k)=ak+b=6k+1.$$
A direct computation shows that in this situation the equality $I_f=(g,g')$ holds no longer.
\end{ex}

\subsection{Cones over hypersurfaces with isolated singularities} 

Consider the polynomial ring $R=\C[x_1,...,x_n]$, a homogeneous polynomial $g \in R$ and
 the  hypersurface $W:g(x_1,...,x_n)=0$  with isolated singularities in $\PP^{n-1}$ with total Tjurina number $\tau(W)$ and stabilization threshold $st(W)$.
Consider the projective cone over $W$, that is the hypersurface $V$ in $\PP^{n}$ defined by
$f=g=0$, where this time $f=g$ is regarded in the polynomial ring $S$. One clearly has
\begin{equation}
\label{tensorp}
M(f)=M(g) \otimes \C[x_0]
\end{equation}
which implies 
$$H(M(f))(k)=\sum_{j=0}^{k}H(M(g))(j),$$
for any positive integer $k$. This yields the following.

\begin{prop}
\label{prop3}
The Hilbert polynomial of the projective cone $V$ over the hypersurface $W$ is given by the formula $P(M(f))(k)=ak+b$ with $a= \tau(W)$ and $$b= \sum_{j=0}^{st(W)-1}H(M(g))(j) -(st(W)-1)\tau(W).$$
Moreover, $st(V)=st(W)-1.$

\end{prop}

\begin{ex}
\label{ex2}
(i) Let $W$ be a curve of degree 4 with 3 nodes. If the nodes are collinear, e.g. 
$$W:g=x(x^3+y^3+z^3)=0$$
then the Hilbert function of $M(g)$ is computed in Example 4.4 in \cite{DStEdin} and we see that 
$st(W)=6$ and $\sum_{j=0}^{st(W)-1}H(M(g))(j)=27$ and hence $b=27-3 \cdot 5=12.$

When the nodes are not collinear, e.g. 
$$W: g= x^2y^2+y^2z^2+x^2z^2-2xyz(x+y+z)-(2xy+3yz+4xz)^2=0$$
then the Hilbert function of $M(g)$ is also computed in Example 4.4 in \cite{DStEdin} and we get  
$st(W)=5$ and $\sum_{j=0}^{st(W)-1}H(M(g))(j)=23$ and hence $b=23-3 \cdot 4=11.$
In both cases the cone $V$ over $W$ has 3 singular lines, each with transversal singularity type $A_1$. In the first case, these 3 lines are in the same plane, determined by the line containing the nodes and the vertex of the cone.

\medskip

\noindent (ii)  Let $W$ be a curve of degree 6 with six cusps $A_2$. If the cusps are on a conic,
then the Hilbert function of $M(g)$ is computed in Example 3.2 in \cite{DStDoc} and we see that 
$st(W)=10$ and $\sum_{j=0}^{st(W)-1}H(M(g))(j)=119$ and hence $b=119-9 \cdot 12=11.$
 If the cusps are not on a conic,
then the Hilbert function of $M(g)$ is also computed in Example 3.2 in \cite{DStDoc} and we see that 
$st(W)=9$, $\sum_{j=0}^{st(W)-1}H(M(g))(j)=105$ and hence $b=105-8 \cdot 12=9.$

\medskip

\noindent (iii) Consider the following two distinct realizations of the  configuration $(9_3)$, see \cite[Example 2.15]{DHA}. The first one is the Pappus line arrangement $(9_3)_1$ given by 
$$W_1: g_1=xyz(x-y)(y-z)(x-y-z)(2x+y+z)(2x+y-z)(-2x+5y-z)=0,$$
with the following Hilbert-Poincar\'e series
$$HP(M(g_1))(t)=1+3t+6t^2+10t^3+15t^4+21t^5+28t^6+36t^7+42t^8+46t^9+48t^{10}+$$
$$+48t^{11}+47t^{12}+45(t^{13}+...$$
It follows that $st(W_1)=13$, $\sum_{j=0}^{st(W)-1}H(M(g_1))(j)=351$ and hence $b=351-12 \cdot 45=-189.$
The second one is the non-Pappus line arrangement  $(9_3)_2$ given by
$$W_2: g_2=xyz(x+y)(x+3z)(y+z)(x+2y+z)(x+2y+3z)(2x+3y+3z)=0.$$
with the following Hilbert-Poincar\'e series
$$HP(M(g_2))(t)=1+3t+6t^2+10t^3+15t^4+21t^5+28t^6+36t^7+42t^8+46t^9+48t^{10}+$$
$$+48t^{11}+46t^{12}+45(t^{13}+...$$
It follows that $st(W_2)=13$, $\sum_{j=0}^{st(W)-1}H(M(g_1))(j)=350$ and hence $b=350-12 \cdot 45=-190.$ Note that $W_1$ and $W_2$ consist both of 9 lines, and have the same number of double and triple points, i.e. 9 double points and 9 triple points.
The fact that the monomial $t^9$ occurs with the same coefficient in both $HP(M(g_1))(t)$ and $HP(M(g_2))(t)$, implies that the two orbits $G \cdot g_1$ and $G \cdot g_2$ , where $G=G\ell_3(\C)$, have the same codimension in $S_9$, namely 46. In particular, none of these two curves $W_1$ and $W_2$ can be regarded as a specialization of the other.

\end{ex}
These examples imply the following surprizing fact.

\begin{cor}
\label{corA}
The Hilbert polynomial of the Milnor algebra $M(f)$ of a projective cone depends on the position of the singularities of the hypersurface $W:g=0$, even in the case when $W$ is a line arrangement.
\end{cor}

\begin{prop}
\label{prop3.2}
The projective cone $V:f=0$ over the hypersurface $W:g=0$ satisfies $N(f)=0$.
\end{prop}

\proof
The equation \eqref{tensorp} shows that $M(f)$ is a free $\C[x_0]$-module, in particular the multiplication by $x_0$ is injective. Since $N(f) \subset \ker \{x_0^N: M(f) \to M(f)\}$ for large $N$, the result follows.
\endproof

\begin{prop}
\label{prop3.3}
Consider the projective cone $V:f=0$ over the hypersurface $W:g=0$. If $K^*$ is a minimal resolution of the graded $R$-module $M(g)$, then $K^*\otimes \C[x_0]$ is a minimal resolution of the graded $S$-module $M(f)$.

\end{prop}

\proof Note that $f_0=0$, so we can chose $f_1,...,f_n$ as a set of generators for $J_f$.
We show only that a syzygy
$a_1f_1+a_2f_2+...+a_nf_n=0$
with coefficients $a_i \in S$ is a linear $\C[x_0]$-combination of syzygies with coefficients in $R$.
To see this, it is enough to write each coefficient $a_i$ as a polynomial in $x_0$, namely
$$a_i=a_{i0}+a_{i0}x_0+...+a_{ik_i}x_0^{k_i}$$
with coefficients $a_{ij} \in R$, and then look at the coefficients of the various powers of $x_0$.
\endproof

\begin{rk}
\label{rkAR}
(i) In terms of the graded modules of all relations, defined by
$$AR(g)=\{(a_1,...,a_n) \in R^n~~:~~a_1g_1+...a_ng_n=0\}$$
and
$$AR(f)=\{(a_0,a_1,...,a_n) \in S^{n+1}~~:~~a_0f_0+a_1f_1+...a_nf_n=0\}$$
one clearly has
$$AR(f)=S\cdot e_0+\C[x_0]AR(g) \subset S \times S^n=S^{n+1}.$$
Here $e_0=(1,0,0,..,0)$ and the components of an element of $AR(g)$ are placed in the second factor $S^n$. 

\noindent (ii) It is clear that conversely $V:f=0$ is a cone over a hypersurface $W:g(x_1,...,x_n)=0$   in $\PP^{n-1}$
if the Jacobial ideal $J_f$ can by generated by $n$ elements, i.e. if the minimal resolution of $M(f)$ ends with the terms
$$S(-(d-1))^n \to S.$$
\end{rk}

The following variant of the cone construction is also useful.
Consider again the polynomial ring $R=\C[x_1,...,x_n]$, a homogeneous polynomial $g \in R$ of degree $d-1$ and
 the  hypersurface $W:g(x_1,...,x_n)=0$  with isolated singularities in $\PP^{n-1}$ and total Tjurina number $\tau(W)$ and stabilization threshold $st(W)$.
Consider  the hypersurface $V$ in $\PP^{n}$ defined by
$$f(x_0,x_1,...,x_n)=x_0g(x_1,...,x_n)=0.$$
 Since $S=R[x_0]$, one clearly has the following decomposition of $M(f)$ as a countable direct sum of obvious $\C$-vector spaces
\begin{equation}
\label{directsum}
M(f)=R/(g) \oplus \oplus_{j \geq 1}M(g)x_0^j.
\end{equation}
This yields the following result for $n=3$ (for $n>3$ the hypersurface $V$ has a singular locus of dimension $>1$, and a similar result can be proved).

\begin{prop}
\label{prop3.5}
For $n=3$, the Hilbert polynomial of the surface $V:  x_0g=0$  in $\PP^3$ is given by the formula $P(M(f))(k)=ak+b$ with $a= \tau(W)+d-1$ and $$b= \sum_{j=0}^{st(W)-1}H(M(g))(j) -st(W)\tau(W) -\frac{(d-4)(d-1)}{2}.$$
Moreover, $st(V)=\max \{d-1,st(W)-1\}.$

\end{prop}

\begin{rk}
\label{rkAR2}
If we look at the syzygies of $f$, we note that all the syzygies of $g$ are still present and we need an additional syzygy of degree 1, namely
$$(d-1)x_0f_0-\sum_{k=1}^{n}x_kf_k=0.$$
This new syzygy is clearly independent of the old ones (i.e. there are no new second order syzygies), and so given the resolution for $M(g)$ it is obvious how to get the resolution for $M(f)$.
\end{rk}

\begin{ex}
\label{ex2.5}
 Consider again the two distinct realizations of the configuration $9_3$ from Example \ref{ex2}. The surface $V_1:f_1=0$ associated to the arrangement $W_1:g_1=0$ has the Hilbert polynomial $P(M(f_1))(k)=54k-261$. Moreover the minimal resolution for $M(g_1)$ is given by
$$0 \to R(-15)^2 \to R(-12)\oplus R(-14)^3 \to R(-8)^3 \to R,$$
and the minimal resolution for $M(f_1)$ is given by
$$0 \to S(-16)^2 \to S(-10) \oplus S(-13)\oplus S(-15)^3 \to S(-9)^4 \to S.$$
The corresponding data for the surface $V_2:f_2=0$ are the following.
$$P(M(f_2))(k)=54k-262,$$
$$0 \to R(-15) \to  R(-13)^3 \to R(-8)^3 \to R,$$
and
$$0 \to S(-16) \to S(-10) \oplus S(-14)^3 \to S(-9)^4 \to S.$$
The interest of considering this construction instead of the cone construction is that the resulting plane arrangements in $\PP^3$ are now essential, i.e. the intersection of all the planes in the arrangement is the empty set.
\end{ex}

\section{Free surfaces in $\PP^3$} 

Recall that the reduced hypersurface $V:f=0$ in $\PP^n$ is said to be free if the module of all relations $AR(f)$) is a free graded  $S$-module. In such a case its rank is $n$ and if
$$r_i=(r_{i0},..., r_{in}) \in AR(f) \subset S^{n+1}$$
for $i=1,...,n$ is a homogenoeous basis of $AR(f)$ with $\deg r_i=d_i$, then the integers $d_i$ are called the {\it exponents} of $V$ (or of $f$).

A key result is the following Saito's criterion:  $V:f=0$ is free if and only if one can find homogeneous elements $r_i \in AR(f)$ for $i=1,...,n$ such that 
\begin{equation}
\label{Saito}
\phi(f)=c \cdot f,
\end{equation}
where  $\phi(f)$ is the determinant  of the $(n+1) \times (n+1)$ matrix $\Phi (f)=(r_{ij})_{0 \leq i,j \leq n}$ and $c$ is a nonzero constant, see for instance \cite[Theorem 8.1]{DHA}. Here the first line in $\Phi(f)$ is
\begin{equation}
\label{Euler}
r_0=(r_{00},..., r_{0n})=(x_0,x_1,...,x_n).
\end{equation}
This vector is not in $AR(f)$, but it corresponds to the Euler derivation. 
If  the condition \eqref{Saito} holds, then the $(r_i)_{i=1,n}$ are a basis for the free $S$-module $AR(f)$ and hence one gets
$$d_1+d_2+...+d_n=d-1.$$
For more on this see \cite{KS}, \cite{Yo}.

We have the following obvious consequence of 
Proposition \ref{prop3.3}.

\begin{cor}
\label{corFC}
 If $W$ is a free divisor in $\PP^{n-1}$, with a minimal resolution for $M(g)$ given by
$$ 0 \to \oplus_{j=1}^{n-1}R(-(d-1)-d_j) \to R(-(d-1))^n \to R,$$
 then $V$ is a free divisor in $\PP^n$ , with a minimal resolution for $M(f)$ given by
$$ 0 \to \oplus_{j=1}^{n-1}S(-(d-1)-d_j) \to S(-(d-1))^n \to S.$$
\end{cor}

In terms of exponents, Corollary \ref{corFC}  says that the exponents of the free divisor $V$ which is a cone over a free divisor $W$ in $\PP^{n-1}$ are given by $d_0=0$, $d_1,...,d_{n-1}$, where $d_1,...,d_{n-1}$ are the exponents of the divisor $W$.

The modified cone construction discussed in Remark \ref{rkAR2} yields the following.

\begin{cor}
\label{corFC2}
 If $W:g=0$ is a free divisor in $\PP^{n-1}$, with a minimal resolution for $M(g)$ given by
$$ 0 \to \oplus_{j=1}^{n-1}R(-(d-2)-d_j) \to R(-(d-2))^n \to R,$$
 then $V: f=x_0g=0$ is a free divisor in $\PP^n$ , with a minimal resolution for $M(f)$ given by
$$ 0 \to S(-d) \oplus \oplus_{j=1}^{n-1}S(-(d-1)-d_j) \to S(-(d-1))^{n+1} \to S.$$
\end{cor}

In terms of exponents, Corollary \ref{corFC2}  says that the exponents of the free divisor $V$ are given by $d_0=1$, $d_1,...,d_{n-1}$, where $d_1,...,d_{n-1}$ are the exponents of the divisor $W$.

In the case of curves, the freeness condition is equivalent to 
$H^0_Q(M(f))=0$, where $H^0_Q$ denotes the 0-th local cohomology module with respect to the ideal  $Q=(x,y,z)$ in $S=\C[x,y,z]$. 
Moreover, in the curve case, the singularities are isolated, since we consider only reduced curves.
For the next result, in the local analytic context,  see also \cite{BP}, section 2, especially Prop. 2.13.

\begin{prop}
\label{propFREE}
The surface $D:f=0$ is free if and only if
$$H^0_Q(M(f))=H^1_Q(M(f))=0,$$
 with $Q=(x,y,z,w)$ the maximal homogeneous ideal in in $S=\C[x,y,z,w]$.
\end{prop}

\proof
By definition, it is clear that $D$ is free if and only if the $S$-module $M(f)$ has projective dimension $\pd M(f)=2$, i.e. $M(f)$ has a minimal free resolution of length 2, see  \cite{Eis}, Cor.1.8. The Auslander-Buchsbaum formula, see \cite{Eis} Thm. A2.15, implies that then $\depth M(f)=\depth S-\pd M(f)=4-2=2.$  Finally \cite{Eis}, Thm. A2.14 tells us that
$$\depth M(f)= \inf \{k ~~|~~H^k_Q(M(f))\ne 0 \}.$$

\endproof

It is known that $N(f)=H^0_m(M(f))=I_f/J_f$, so the corresponding vanishing is easy to check using a computer software. The second vanishing is more subtle. Using the definition of the defect and the formula \eqref{Hfunction2}, we get the following.

\begin{cor}
\label{cor4}
 If a surface $D:f=0$  in $\PP^3$ is free, then $D$ has a 1-dimensional singular locus $\Sigma$.
Conversely, a surface $D:f=0$  in $\PP^3$ with a 1-dimensional singular locus $\Sigma$ is free if and only if $N(f)=I_f/J_f=0$ and one of the following equivalent conditions holds.

\noindent (i) $P(f)_k=\defect_k\Sigma=0$ for any $k$;

\noindent (ii) The evaluation morphism $ev: S_k/I_{f,k} \to H^0(\Sigma, \OO_{\Sigma}(k))$ is surjective for any $k$; in particular, the support $|\Sigma|$ of the singular locus has to be connected.

\noindent (iii) $H(M(f))(k)=\dim H^0(\Sigma, \OO_{\Sigma}(k))$ for any $k$.

\end{cor}

For the fact that a free divisor $D$ has a 1-dimensional singular locus, see \cite{BP}, section 2, especially Prop. 2.13 or Corollary \ref{corFREE} (ii) below.

The property (ii) is known as the $k$-normality of the projective (in general non-reduced) curve $\Sigma$.
One knows that if $\Sigma$ is a complete intersection, i.e. if the ideal $I_f$ is generated by two polynomials, then the condition (ii) above is fullfilled, see \cite{FAC}, Proposition 5, p. 273-274. 
However, this yields free divisors practically never, as the following shows.

\begin{ex}
\label{exNONF}
A surface $V=D \cup D':f=gg'=0$ as in subsection \ref{normalcross} with $(e,e') \ne (1,1)$ is not free (since clearly  $N(f) \ne 0$ by looking at the degree of the generators), although by Proposition \ref{prop2} the condition (ii) is fulfilled. This is similar to the fact that a nodal plane curve of degree $d\geq 4$ cannot be free, see \cite{DS14}, Example 4.1 (i).
\end{ex}

We recall now the definition of some invariants associated with a projective surface in $\PP^3$, see \cite{DStEdin}.

\begin{definition}
\label{def0}

For a reduced surface $D:f=0$ of degree $d$ in $\PP^3$, two invariants are defined as follows.

\noindent (i) the {\it coincidence threshold} 
$$ct(f)=\max \{q:\dim M(f)_k=\dim M(f_s)_k \text{ for all } k \leq q\},$$
with $f_s$  a homogeneous polynomial in $S$ of degree $d$ such that $D_s:f_s=0$ is a smooth surface in $\PP^3$.

\noindent (ii) the {\it minimal degree of a  nontrivial (or essential) syzygy} 
$$mdr(f)=\min \{q~~:~~ H^3(K^*(f))_{q+3}\ne 0\},$$
where $K^*(f)$ is the Koszul complex of $f_x,f_y,f_z,f_w$ with the natural grading.
\end{definition}
Note that one has  for $j<d-1$ the following equality
\begin{equation} 
\label{ar=er}
AR(f)_j=H^3(K^*(f))_{j+3}.
\end{equation} 
Moreover, $mdr(f)=0$ if and only if $f$ is independent of $x_0$ after a linear change of coordinates, i.e. $D:f=0$ is a cone over a curve in $\PP^{2}$.
It is known that one also has
\begin{equation} 
\label{REL}
ct(f)=mdr(f)+d-2.
\end{equation} 

\begin{thm} \label{thmFREE}
 Suppose the surface $D:f=0$ is not a cone and it is free with exponents $1\leq d_1 \leq d_2 \leq d_3$, i.e. the minimal resolution of the Milnor algebra has the form
$$0 \to \oplus_{j=1}^{3}S(-(d-1)-d_j) \to S(-(d-1))^4 \to S.$$
Then one has the following.

\noindent (i) $d_1+d_2 +d_3= d-1$, $mdr(f)=d_1$, $ct(f)=d_1+d-2$ and $st(f) =d+d_3-4$.

\noindent (ii) The coefficients of the Hilbert polynomial $P(M(f))(k)=ak+b$ are given by
$$a=s_1^2-s_2$$
and
$$b=2a-s_1^3+\frac{3}{2}s_1s_2-\frac{1}{2}s_3,$$
where $s_1=\sum_{j=1}^{3}d_j$, $s_2=\sum_{i <j}d_id_j$ and $s_3=d_1d_2d_3$ are the elementary symmetric functions in the exponents.
\end{thm}

\proof
The given minimal resolution of the Milnor algebra implies the equality
$$\dim M(f)_k={k+3 \choose 3}-4{k+4-d \choose 3}+\sum_{j=1}^{3}{k+4-d-d_j \choose 3}$$
for all $k \geq d+d_3-4$. If we expand the right hand side we get a polynomial $Q(k)$ in $k$ of degree 2.
Since $P(M(f))(k)$ is a polynomial of degree  at most 1 (as $D$ is reduced,  the singular locus of $D$ is at most 1-dimensional) it follows that the leading coefficient of $Q(k)$ must vanish and this implies the equality  $s_1= d-1$. The other equalities in the claim (i) follow by definition.

To get the formulas in the claim (ii), we identify the coefficients of the polynomials $P(M(f))(k)$ and $Q(k)$ and write the result in terms of the exponents $d_i$, using the equality  $s_1= d-1$. 

It is easy to show that the expression for $b$ is an integer  for any exponents $d_i$'s, so this result imposes no condition on the exponents except $s_1= d-1$. 

\endproof

\begin{cor}
\label{corFREE}
(i) The exponents of a free divisor  $D:f=0$ in $\PP^3$ which is not a cone determine the degree $d$ of $D$ and the Hilbert polynomial $P(M(f))$.
Conversely, the degree $d$ of a such a free divisor $D$ and the Hilbert polynomial $P(M(f))$ determine
the exponents.

(ii) The degree of the singular locus $\Sigma$ of a free divisor which not a cone is given by
$$\deg \Sigma=a=\sum_{j=1}^{3}d_j^2+\sum_{i <j}d_id_j.$$
Moreover, one has $ \deg \Sigma \geq 6$, in particular $\dim \Sigma=1$.

(iii) If $D:f=0$ is a free divisor, then $$\dim H^1(\Sigma, \OO_{\Sigma}(k))=H(M(f))(k)-P(M(f))(k)$$
for any $k\in \Z$.
\end{cor}

\proof 

The first two claims are clear by Theorem \ref{thmFREE}. The last follows from \eqref{Hfunction2}.

\endproof

The reader may state the corresponding properties (i) and (ii) in the case of a free surface which is a cone, whose resolution is a special case of Corollary \ref{corFC} for $n=3$. The case when $W$ is a cone itself, i.e. a union of lines passing through one point, has to be treated separatedly.

\begin{ex}
\label{exFREE}
(i) The plane arrangement $D=xyzw=0$ is free with exponents $d_1=d_2=d_3=1$ and the corresponding Hilbert polynomial is $P(M(f))(k)=6k-2$. The corresponding Hilbert function is $H(M(f))(0)=1$ and $H(M(f))(k)=P(M(f))(k)$ for $k>0$. For more on free hyperplane arrangements see \cite{OT}, \cite{Yo}.

\noindent (ii) The plane arrangement $D=(x^3+y^3)(z^3+w^3)=0$ is free with exponents $d_1=1$, $d_2=d_3=2$ and the corresponding Hilbert polynomial is $P(M(f))(k)=17k-33$.
The corresponding Hilbert function is $H(M(f))(0)=1$,  $H(M(f))(1)=4$, $H(M(f))(2)=10$, $H(M(f))(3)=20$ and $H(M(f))(k)=P(M(f))(k)$ for $k>3$.

\end{ex}

As in the case of free plane curves, to find irreducible examples is harder. Here are some  irreducible free surfaces.

\begin{ex}
\label{exFREE2}
(i) Consider the sequence of surfaces 
$$D_d:f_d=x^{d-1}z+y^d+x^{d-2}yw+x^4y^{d-4}$$
for $d \geq 4$.
The surface corresponding to $d=7$ is free with exponents $d_1=1$, $d_2=2$, $d_3=3$,
and the corresponding Hilbert polynomial is $P(M(f_7))(k)=25k-70$. The Hilbert series of $M(f_7)$ is
$$HP(M(f_7;t)=1+4t+10t^2+20t^3+35t^4+56t^5+...$$
with $st(f_7)=6$.
The surface corresponding to $d=8$ is also free with exponents $d_1=1$, $d_2=3$, $d_3=3$,
and the corresponding Hilbert polynomial is $P(M(f_8))(k)=34k-122$. The Hilbert series of $M(f_8)$ is
$$HP(M(f_8;t)=1+4t+10t^2+20t^3+35t^4+56t^5+84t^6+...$$
with $st(f_8)=7$. The remaining surfaces $D_d$  are discussed below in Proposition \ref{propDd}.

(ii) Consider the sequence of surfaces 
$$D'_d:f'_d=x^{d-1}z+y^d+x^{d-2}yw+x^{d-5}y^5$$
for $d \geq 6$. These surfaces are free for $d\geq 10$, with exponents $d_1=1$, $d_2=4$, $d_3=d-6$, see Proposition \ref{propD'd} below.
The surfaces $D'_k$ for $6\leq k \leq 9$ are discussed below in Example \ref{exNF} (vi).

\end{ex}

\begin{prop}
\label{propD'd}
The surface
$$D'_d:f'_d=x^{d-1}z+y^d+x^{d-2}yw+x^{d-5}y^5$$
 is free with exponents $(1,4,d-6)$ for any $d \geq 10$.

\end{prop}

\proof
Looking at the results produced by Singular for various $d$, with $10 \leq d \leq 20$, 
we find out the following expressions for the generators of the first syzygies.
$$r_1=( 0, 0, -y,  x),$$
$$r_2=(x^4, 0, -(d-1)x^3z-(d-2)x^2yw, -(d-5)y^4),$$
$$r_3=( r_{3x},r_{3y},r_{3z},r_{3w})$$
with
$$r_{3x}=-d(d-5)y^{d-6} + d(d-2)x^3y^{d-10}w,$$
$$r_{3y}=
(d-5)^2x^{d-6}, $$ 
$$r_{3z}=
-5(d-5)^2x^{d-10}y^4 - (d-5)^2x^{d-7}w - d(d-1)(d-2)x^2y^{d-10}zw -d(d-2)^2xy^{d-9}w^2,$$
and 
$$r_{3w}=
d(d-1)(d-5)y^{d-7}z .$$
To prove formally that $D'_d$ is a free divisor, we may now either use Saito's criterion \eqref{Saito}, or note that the following version of Lemma 1.1 in \cite{ST} holds with the same proof.

\endproof

\begin{lem}
\label{lemST}
The surface $D:f=0$ in $\PP^3$ is free if and only if there exist three distinct minimal generating syzygies of degrees $e_1$, $e_2$ and $e_3$ such that $e_1+e_2+e_3 \leq d-1$.
\end{lem}

Next we discuss some discriminants, objects that played a key role in the introduction by K. Saito of the notion of free divisors in \cite{KS}. 
Let $\Delta_n \in \C[a_0,...,a_n]$ be the discriminant of the
general degree $n$ binary form
$$a_0x^n+a_1x^{n-1}y+...+a_{n-1}xy^{n-1}+a_ny^n.$$
One has the following result, see \cite{BEG} and especially \cite{DCA}.

\begin{prop}
\label{propDn}

For $n\geq 3$, the hypersurface $\D_n:\Delta_n=0$ in $\PP^n$ has degree $d=2n-2$, is  irreducible and free  with exponents $d_1=d_2=d_3=1$, $d_4=...=d_n=2$.

\end{prop}
In particular,  $\D_3$ is a linear free divisor,  for more on this see \cite{BM} and  \cite{GS}, Example 1.6.
In this series of free divisors, only the first two play a role in this note, and we give below their known equations and some additional properties.

\begin{prop}
\label{propD3}
The discriminant $\Delta_3$ of the binary form 
$$P=ax^3+bx^2y+cxy^2+dy^3$$
is given by 
$$\Delta_3=b^2c^2 - 4ac^3 - 4b^3d + 18abcd - 27a^2d^2  .$$
Moreover,  the surface $\D_3$ is rational, homeomorphic to $\PP^1 \times \PP^1$, and has only points of multiplicity $\leq 2$.

\end{prop}

\proof

Note that a binary form $P$ is in $\D_3$ if and only if it has a double factor, hence if and only if it can be written as $L_1^2L_2$, with $L_1,L_2$ linear forms in $x$ and $y$.
It follows that the morphism $\PP^1 \times \PP^1 \to \D_3$ sending $(L_1,L_2)$ to $L_1^2L_2$
is a bijection. It follows that $\D_3$ is rational, homeomorphic to $\PP^1 \times \PP^1$.

The points of highest multiplicity on $\D_3$ corresponds to a perfect cube, e.g. the point
corresponding to $y^3$ (all such points are in a $Sl_2(\C)$ orbit, so they have the same multiplicity). The equation of $\D_3$ implies that the multiplicity of $y^3$ is 2.
Note that all the other examples of irreducible degree $d$ free divisors given above have points of multiplicity $d-1$, and hence their rationality is obvious.

\endproof

\begin{prop}
\label{propD4}
The discriminant $\Delta_4$ of the binary form 
$$P=ax^4+bx^3y+cx^2y^2+dxy^3+ey^4$$
is given by 
$$\Delta_4= b^2c^2d^2 - 4ac^3d^2 - 4b^3d^3 + 18abcd^3 - 27a^2d^4 - 4b^2c^3e + 16ac^4e + 18b^3cde - $$
$$80abc^2de - 6ab^2d^2e + 
144a^2cd^2e - 27b^4e^2 + 144ab^2ce^2 - 128a^2c^2e^2 - 192a^2bde^2 + 256a^3e^3.$$
Moreover, the singular locus of 3-fold $\D_4:\Delta_4=0$  has dimension 2 and the corresponding Hilbert polynomial is 
$$H(M(f))(k) = 8k^2 - 21k + 26.$$
The multiplicity of points on $\D_4$ is $\leq 3$.

\end{prop}

\proof
Exactly as the previous proof.
\endproof
The 3-fold $\D_4$ and its defining equation play a key role in Proposition \ref{propD4'} below, which explain their inclusion in our discussion, though $\D_4$ has a 2-dimensional singular locus.

\section{Nearly free surfaces in $\PP^3$} 

First we recall the definition of nearly free curves.

\begin{definition} 
\label{def2}

The curve $C:f=0$ is a { nearly free divisor} 
if the following  equivalent conditions hold.

\begin{enumerate}

\item[(i)] $N(f) \ne 0$ and $n(f)_k \leq 1$ for any $k$.

\item[(ii)] The Milnor algebra $M(f)$ has a minimal resolution of the form
$$0 \to S(-d-d_2) \to S(-d-d_1+1) \oplus S^2(-d-d_2+1) \to S^3(-d+1) \xrightarrow{(f_0,f_1,f_2)}  S$$
for some integers $1 \leq d_1 \leq d_2$, called the exponents of $C$.

\end{enumerate}

\end{definition}

\begin{rk}
\label{rkB+}
The nearly free curves were introduced in \cite{DStNF}, where we conjectured that any rational cuspidal curve in $\PP^2$ is either free, or nearly free, and we proved this conjecture for all the curves of even degree, and for many odd degree curves. In the case of plane curves, the definition of nearly freeness was dictated by the aim of constructing a minimal class of curves, such that the above conjecture holds. Moreover, the free and nearly free curves share many interesting properties, e.g. they maximize the total Tjurina number, when we fix the degree $d$ and the minimal degree of a Jacobian syzygy, see \cite{duPCTC} for the free curves and \cite{DmaxC} for the nearly free curves.
A secondary conjecture in \cite{DStNF}, based on the known examples of irreducible free curves at that moment, was that a free curve is necessarily rational and not far from a cuspidal curve. This conjecture was disproved in the paper \cite{B+}.
\end{rk}

After this motivation of our interest in nearly free curves, we give  the definition of a nearly free divisor in $\PP^3$. 

\begin{definition}
\label{NFdef}
The reduced surface  $D:f=0$ in $\PP^3$ is nearly free if either

\begin{enumerate}

\item[(i)] $ mdr(f)=0$, i.e. $D$ is a cone over a nearly free curve $C$ in $\PP^2$, and then the minimal resolution of the Milnor algebra $M(f)$ has the form
\begin{equation}
\label{rnfc}
 0 \to S(-d_2-d) \to S(-d_1-(d-1))  \oplus S(-d_2-(d-1))^2 
\to S(-(d-1))^3 \to S
\end{equation}
for some integers $1 \leq d_1 \leq d_2 $, or

\item[(ii)]  $mdr(f)>0$, i.e. $D$ is not a cone over a plane curve, and the minimal resolution of the Milnor algebra $M(f)$ has the form
\begin{equation}
\label{rnfs}
 0 \to S(-d_3-d) \to S(-d_1-(d-1)) \oplus S(-d_2-(d-1)) \oplus S(-d_3-(d-1))^2 \to  
\end{equation}
$$\to S(-(d-1))^4 \to S$$
for some integers $1 \leq d_1 \leq d_2 \leq d_3$.
\end{enumerate}
 
\end{definition}

In down-to-earth terms, this definition in case (ii) says that the module $AR(f)$ is not free of rank 3, but has 4 generators $r_1$, $r_2$, $r_3$ and $r_4$ of degree respectively $d_1$, $d_2$, $d_3$ and $d_3$ and the second order syzygies are spanned by a unique relation
\begin{equation}
\label{SOS}
R: a_1r_1+a_2r_2+a_3r_3+a_4r_4=0.
\end{equation}
Here $a_1,a_2,a_3,a_4$ are homogeneous polynomials in the graded polynomial ring $S$ of degrees $d_3-d_1+1, d_3-d_2+1,1,1$ respectively.

One has the following version of Saito's criterion \eqref{Saito} in the case of a nearly free divisor.
 Recall that $r_0$ corresponds to the Euler vector field \eqref{Euler}.

\begin{thm}
\label{thmSaito}
Let $D:f=0$ be a nearly free surface which is not a cone and let $r_i\in AR(f)$ and $a_j\in S$ be as above.
For $i=1,2,3,4$ denote by $\phi_i$ the determinant of the $4 \times 4$ matrix whose lines are in order given by the components of the first order syzygies $r_0$,..., $r_{i-1}$, $r_{i+1}$,...,$r_4$.
Then there is a nonzero constant $c$ such that 
$$\phi_i=(-1)^i ca_if$$ 
for $i=1,2,3,4$. Conversely, let $D:f=0$ be a reduced surface which is not a cone and assume that

\noindent (i) The Jacobian ideal $J_f$ is saturated.

\noindent (ii) The module $AR(f)$ has a set of minimal homogeneous generators $r_1,r_2,r_3,r_4$ of degrees $d_1 \leq d_2 \leq d_3=d_4$ such that $d_1+d_2+d_3=d$, the degree of $D$.

Then $D:f=0$ is a nearly free surface and the coefficients $a_i$ of the second order syzygy \eqref{SOS} can be taken to be $a_i=(-1)^i\phi_i/f$ for $i=1,...,4$.
\end{thm}

\proof
Consider the matrix $A$ with 5 rows and 4 columns obtained by putting in order the rows corresponding to $r_i$ for $i=0,1,...,4$. Construct a square matrix $B$, by adding a fifth column to $A$, consisting in order of $(0,r_{1,x},r_{2,x},r_{3,x},r_{4,x})$, where of course $r_{1,x}$ denotes the first coordinate of $r_1$ and so on. It is clear that $\det B=0$, since the last 4 rows in $B$ are linearly dependent over the field of fractions $K(S)$ of $S$, due to the relation \eqref{SOS}.
Expanding this determinant over the last column, we get
$$r_{1,x}\phi_1-r_{2,x}\phi_2+r_{3,x}\phi_3-r_{4,x}\phi_4=0.$$
Note that all the determinants $\phi_i$ vanish on the smooth part of $D$ (since the rows of the matrix correspond to tangent vectors to $D$), and hence one has $\phi_i=f \phi_i'$, for some polynomials $\phi_i'$. Hence we have
$$r_{1,x}\phi_1'-r_{2,x}\phi_2'+r_{3,x}\phi_3'-r_{4,x}\phi_4'=0.$$
Note that the last column in $B$ could as well be taken to be $(0,r_{1,y},r_{2,y},r_{3,y},r_{4,y})$, or
$(0,r_{1,z},r_{2,z},r_{3,z},r_{4,z})$, or $(0,r_{1,w},r_{2,w},r_{3,w},r_{4,w})$, and the result would have been the same. It follows that we have in fact
\begin{equation}
\label{SOS'}
\phi_1'r_1-\phi_2'r_2+\phi_3'r_3-\phi_4'r_4=0.
\end{equation}
Since clearly the homogeneous polynomials $\phi_i'$ and $a_i$ have the same degree, it follows that they differ by a multiplicative constant in view of the unicity of the relation \eqref{SOS}.
The relations $r_0$, $r_1$, $r_2$ and $r_3$ are linearly independent over $K(S)$, and this implies $c\ne 0$.

To prove the converse claim, note that the assumption (ii) gives an exact sequence
$$S(-d_1-d+1) \oplus S(-d_2-d+1) \oplus S(-d_3-d+1)^2  
\to S(-d+1)^4 \to S \to M(f) \to 0.$$
The second order syzygy \eqref{SOS'} shows that we can extend the above exact sequence to a complex similar to \eqref{rnfs}.  Note that this syzygy, having two components of degree one, cannot be the multiple of another second order syzygy of lower degree.
It remains to explain why the obtained complex is exact in the second term, i.e. it is a minimal resolution for $M(f)$.
If exactness fails there, it means that we need more second order syzygies. But this implies that the minimal resolution of $M(f)$ has length 4. Arguing as in the proof of Proposition \ref{propFREE}, this leads to a contradiction of the assumption (i). This completes the proof of this theorem.
\endproof

\begin{rk}
\label{rkkeyNF}
(i) The fact that $r_1,r_2,r_3,r_4$ form a system of generators for the module $AR(f)$ is difficult to check in practice. It follows from Proposition 5.1 in \cite{Dmax} that this condition can be replaced when
 $f$ is a tame polynomial in the sense of \cite{Dmax}, by the condition that the two linear forms $a_i=(-1)^i\phi_i/f$ for $i=3,4$ are linearly independent. This condition is much simpler to check in practice. An example of this approach is given by the proof of Proposition \ref{propDd} below.
 
\noindent(ii) The syzygies $r_j$ for $j=0, ..., 4$ considered in Theorem \ref{thmSaito}, regarded as derivations of the polynomial ring $S$ in the obvious way, generate the module $D^1(f)$ of such derivations preserving the ideal $(f)$. The relation \eqref{SOS} implies that the determinant constructed using  the components of the  syzygies $r_1$,$r_2$, $r_3$ and $r_4$ is zero. It follows that the ideal $min_4(A)$ spanned by all the possible 4-minors in the $5 \times 4$ matrix $A$  considered in the proof of Theorem \ref{thmSaito} satisfies the equality
$$min_4(A)=(a_1f,a_2f,a_3f,a_4f).$$
One can prove, exactly as in \cite[Theorem 2.11]{TW}, that one has
$$ \sqrt {min_4(A)}=(f).$$
This equality is clearly equivalent to 
 $$ (f) \subset \sqrt {(a_1,a_2,a_3,a_4)}.$$
 For explicit computations of the ideal $I_A=(a_1,a_2,a_3,a_4) \subset S$ see Example \ref{exNF4}.
 \end{rk}

\begin{ex}
\label{exNF}
(i) The surface $D:f=x^2z+y^2w=0$ has the following minimal resolution for $M(f)$
$$0 \to S(-6) \to S(-5)^4 \to S(-3)^2\oplus S(-4)^4 \to S(-2)^4 \to S$$
which is too long, hence this surface is not nearly free.

\noindent (ii) The surface $D:f=x^2z+y^3+xyw=0$ has the following minimal resolution for $M(f)$
$$0  \to S(-5) \to S(-3)^3\oplus S(-4) \to S(-2)^4 \to S$$
which has the right length, but is not of the type in the definition. Hence this surface is not nearly free.

\noindent (iii)  The surface $D:f=x^4-xyw^2+zw^3=0$ has the following minimal resolution for $M(f)$
$$0  \to S(-6) \to S(-4)^2\oplus S(-5)^2 \to S(-3)^4 \to S$$
which is of the type in the definition, with $d_1=d_2=1$ and $d_3=2$. Hence this surface is  nearly free.

\noindent (iv) The surfaces $D'_k$ for $6\leq k \leq 9$  from sequence of surfaces introduced in Example \ref{exFREE2} (ii) are all nearly free with the following exponents.

$\bullet$ for $k=6$ the exponents are $(d_1,d_2,d_3)=(1,2,3)$.

$\bullet$ for $k=7$ the exponents are $(d_1,d_2,d_3)=(1,3,3)$.

$\bullet$ for $k=8$ the exponents are $(d_1,d_2,d_3)=(1,3,4)$.

$\bullet$ for $k=9$ the exponents are $(d_1,d_2,d_3)=(1,4,4)$.

\end{ex}

 We construct now some countable families of  nearly free surfaces. First we consider the family of surfaces $D_d$ introduced in Example \ref{exFREE2} (i).

\begin{prop}
\label{propDd}

The surface $D_d:f=x^{d-1}z+y^d+x^{d-2}yw+x^4y^{d-4}=0$  is nearly free with exponents

\begin{enumerate}

\item  $(1,1,2)$ for $d=4$;

\item  $(1,1,3)$ for $d=5$;

\item  $(1,1,4)$ for $d=6$;

\item  $(1,4,4)$ for $d=9$;

\item  $(1,4,5)$ for $d=10$;

\item $(1,5,d-6)$ for $d\geq 11$.

\end{enumerate}

\end{prop}

\proof

We give the proof in the case $d\geq 11$, the other cases can be settled using a computer algebra system.
One can easily find the first syzygy
$r_1=( 0, 0, -y, x )$, while some work shows that the second syzygy $r_2$ is given by
$$( (d-4)x^5+dxy^4  , -4x^4y  ,  -(d-1)z[(d-4)x^4 +dy^4] - (d^2-6d+4)x^3yw ,  -d(d-2)y^4w   ).$$
The third syzygy $r_3$ is given by 
$$( dx^{d-9}y^3,  -4x^{d-6}, 4x^{d-7}w,  4(d-4)y^{d-6}-d(d-1)x^{d-9}y^2z-d(d-2)x^{d-10}y^3w).$$
Finally, the fourth syzygy $r_4=( r_{40}, r_{41}, r_{42}, r_{43})$ is given by 
$$ r_{40}= d^2(d-2)x^{d-10}y^3w+4(d-4)y^{d-10}[dy^4+(d-4)x^4],     $$
$$ r_{41}=  -[4d(d-2)x^{d-7}w+16(d-4)x^3y^{d-9}],    $$
$$ r_{42}= 4d(d-2)x^{d-8}w^2-4(d-4)x^2y^{d-10}[(d^2-6d+4)yw+(d^2-5d+4)xz],     $$
$$ r_{43}=-[4d(d-1)(d-4)y^{d-7}z+d^2(d-2)(d-1)x^{d-10}y^2zw +d^2(d-2)^2x^{d-11}y^3w^2].      $$
The $2 \times 2$ minors of the matrix formed by the componenets of $r_1$ and $r_2$ have no common divisor in the polynomial ring $S$, which implies that $f$ is a tame polynomial in the sense of \cite{Dmax}, see Lemma 2.2 (ii) in \cite{Dmax}. To apply Remark \ref{rkkeyNF}, one computes 
$a_3= - 4d(d-4)x$ and 
$a_4=-4d^2(d^2-6d+8)w$. This shows that $a_3$ and $a_4$ are linearly independent, which by 
Remark \ref{rkkeyNF} completes the proof that the surface $D_d$ is nearly free for $d\geq 11$.

\endproof

\begin{rk}
\label{rktame}
The class of tame polynomials used in the above proof has also the following property: within this class, and fixing $d=\deg(f)$ and $d_1=mdr(f)$, a surface $D:f=0$ is free if and only if the degree of the associated singular scheme $\Sigma(f)$ is maximal. Hence in this class of surfaces, the analog of the characterization of free curves given in 
\cite{duPCTC, DmaxC} holds, see \cite{Dmax}. Moreover, if $D:f=0$ is free, then $f$ is tame, but it is an open question if the defining equation of a nearly free surface is tame, even if this happens in all examples tested so far.

\end{rk}

\begin{prop}
\label{propD''d}

The surface $D''_d:f=x^{d-1}z+y^d+x^{d-2}yw=0$ is nearly free with exponents $(1,1,d-2)$ for  any $d \geq 4$.

\end{prop}

\proof

 We use the generators of the first syzygies, namely
$$r_1=( x, 0, -(d-1)z, -(d-2)w ) \ \ \ \
r_2=( 0, 0, -y, x ),$$
$$r_3=( dy^{d-2}, -(d-2)x^{d-3}w, (d-2)x^{d-4}w^2, -d(d-1)y^{d-3}z)$$
and
$$r_4=( 0, x^{d-2}, -x^{d-3}w, -dy^{d-2})$$
and Saito's criterion in Theorem \ref{thmSaito} in the form described in Remark \ref{rkkeyNF}.

\endproof

\begin{rk}
\label{rkrepetitions}

We remark that there are some overlaps in our three families $D_d$, $D'_d$ and $D''_d$: up to projective equivalence one has
$D_5=D'_5=D''_5$, $D_6=D''_6$ and $D_9=D'_9$. One can show that these are the only overlaps
(the details will be given elsewhere). Moreover, these surfaces are all homaloidal, i.e. their gradients give rise to birational automorphisms of $\PP^3$, see \cite{DStMM}.

\end{rk}

Now we construct a nearly free divisor having $d_1>1$.
\begin{prop}
\label{propD4'}
Let $\D'_4$ be the surface obtained as the hyperplane section of the free 3-fold $\D_4$ from Proposition \ref{propD4} by the hyperplane $a-c=0$, i.e. the surface given by
$$f=-4a^4d^2 + a^2b^2d^2 + 18a^2bd^3 - 4b^3d^3 - 27a^2d^4 + 16a^5e - 4a^3b^2e - 80a^3bde + 18ab^3de + $$
$$144a^3d^2e - 6ab^2d^2e - 
128a^4e^2 + 144a^2b^2e^2 - 27b^4e^2 - 192a^2bde^2 + 256a^3e^3=0$$
in $\PP^3$ with coordinates $a,b,d,e$.
 Then $\D'_4$ is a nearly free divisor with the exponents $d_1=d_2=d_3=2$. Moreover $\D'_4$ is a rational surface.
\end{prop}

\proof

The only claim which (maybe) is not obtainable by a direct computation is the rationality claim.
Note that the Zariski open set $U_4=\{a=1\} \cap \D_4$ can by identified with the family of polynomials 
$$(x-\alpha)^2(x^2-sx+p).$$
It follows that the map $\phi: \C^3 \to U_4$, $(\alpha,s,p) \mapsto (x-\alpha)^2(x^2-sx+p)$
is an isomorphism above the open set $U_4 \setminus \D_{4,sing}$, since the polynomials with two double roots are clearly in the singular part $\D_{4,sing}$ of $\D_4$.

Consider now the intersection $U_4'=\{a=1\} \cap \D_4'$. It consists of polynomials
$(x-\alpha)^2(x^2-sx+p)$ as above satisfying the extra condition $c=1$, namely
$$p+2\alpha s+\alpha^2=1.$$
It follows that the map $\phi': \C^2 \to U_4'$, 
$$(\alpha,s) \mapsto (x-\alpha)^2(x^2-sx+1-(2\alpha s+\alpha^2))$$
is an isomorphism above the open set $U_4' \setminus \D_{4,sing}'$.
\endproof

\medskip

Exactly the same proof as for Theorem \ref{thmFREE}, but based now on the formula
$$\dim M(f)_k={k+3 \choose 3}-4{k+4-d \choose 3}+\sum_{j=1}^{4}{k+4-d-d_j \choose 3}-{k+3-d-d_3 \choose 3},$$
where $d_4=d_3$ and $k$ is large enough, yields the following result.

\begin{thm} \label{thmNFREE}
 Suppose the surface $D:f=0$ is nearly free with exponents $1\leq d_1 \leq d_2 \leq d_3$, i.e. the minimal resolution of the Milnor algebra has the form given in Definition \ref{NFdef}.
Then one has the following.

\noindent (i) $d_1+d_2 +d_3= d$, $mdr(f)=d_1$, $ct(f)=d_1+d-2$ and $st(f) =d+d_3-3$.

\noindent (ii) The coefficients of the Hilbert polynomial $P(M(f))(k)=ak+b$ are given as follows.
 Define $d'_1=d_1$, $d'_2=d_2$, $d'_3=d_3-1$ and let the integers $a'$ and $b'$ be computed using the formulas in Theorem \ref{thmFREE} (ii), i.e. as if $a',b'$ were the coefficients of the Hilbert polynomial corresponding to a free surface $D'$ with exponents $d'_1,d'_2,d'_3$. Then one has the formulas 
$$a=a'-1 \text{   and   } b=b'+d+d_3-3.$$

\end{thm}

\begin{rk}
\label{rkNF1}

 The invariant $st(f)$ is related to the Castelnuovo-Mumford regularity $\reg M(f)$. By definition, see \cite{Eis}, p. 55, if $D:f=0$ is either free or nearly free, then $\reg M(f)=d+d_3-3$.
Moreover, Theorem 4.2 (2) in \cite{Eis} implies that $st(f) \leq \reg M(f)-\delta+3$, where $\delta$ is the projective dimension of $M(f)$. It follows that $st(f) \leq d+d_3-4$ for a free surface, and $st(f) \leq d+d_3-3$ for a nearly free surface. By our results above, both of these inequalities are in fact equalities. In the case of a free surface, the equality follows also from
Theorem 4.2 (3) in \cite{Eis} since $M(f)$ is a Cohen-Macaulay module in this case, by Theorem A2.14 (4) in \cite{Eis}.
For a nearly free surface,  the last quoted result implies that $M(f)$ is not a Cohen-Macaulay module.

\end{rk}

\begin{cor}
\label{corNF1}
The singular locus $\Sigma$ of a nearly free surface $D$ is 1-dimensional, and of degree at least $d_3$ if $D$ is not a cone.
\end{cor}

\proof

One has the obvious inequalities
$$\deg \Sigma=a=\sum_{j=1}^{3}d_j^2+\sum_{i <j}d_id_j-d -d_3 \geq \sum_{j=1}^{3}d_j+ (d_1+d_2)d_3-d-d_3\geq d_3 >0.$$

\endproof

\begin{ex}
\label{exNF2} 

The computation using Singular gives the same results as  Theorem \ref{thmNFREE} for the nearly free surfaces from Example \ref{exNF}. Here are two cases.

\noindent (i)  The surface $D:f=x^4-xyw^2+zw^3=0$ is nearly free with exponents $d_1=d_2=1$ and $d_3=2$ and the corresponding Hilbert polynomial is $P(M(f))(k)=5k+1$. Moreover $st(f)=3$.

\noindent (ii)  For series of nearly free surfaces $D''_d:f=x^{d-1}z+y^d+x^{d-2}yw=0$ with $d \geq 4$ , the corresponding Hilbert polynomial is $P(M(f))(k)=ak+b,$
where $a=d^2-4d+5$ and $b=-(d^3-8d^2+20d-17)$. Moreover $st(f)=2d-5$.
\end{ex}

Finally we investigate the local cohomology of the Milnor algebra of a nearly free divisor.
\begin{thm}
\label{thmNF}
Let $D:f=0$ be surface  in $\PP^3$ which is not a cone and it is  nearly free. Then the 
following hold.

\noindent (i) $H^0_Q(M(f))=0$, i.e. the Jacobian ideal $J_f$ is saturated.

\noindent (ii) $H^1_Q(M(f))_k = 0$ for any  integer $k > d+d_3-4$ and the following are equivalent

\begin{enumerate}

\item $H^1_Q(M(f))$ is a finite dimensional $\C$-vector space;

\item the coefficients $a_1,a_2,a_3,a_4$ in the relation \eqref{SOS}  form a regular sequence in $S$, i.e. they  define a 0-dimensional complete intersection.

\end{enumerate}

\end{thm}

\proof

The first claim is clear by the length of the resolution \eqref{rnfs} given in Definition \ref{NFdef}, see also the proof of Proposition \ref{propFREE}.
Consider now the dual $L^*$ of this resolution, twist it by $(-4)$ and look at the last morphism
$$\delta_3: L^2(-4)=S(d+d_1-5) \oplus S(d+d_2-5) \oplus S(d+d_3-5)^2 \to L^3(-4)=S(d+d_3-4).$$
This morphism has the form
$$(u_1,u_2,u_3,u_4) \mapsto a_1u_1+a_2u_2+a_3u_3+a_4u_4$$
where $a_1,a_2,a_3,a_4$ are the homogeneous polynomials from \eqref{SOS}.
It is known, see \cite{Eis}, that the dual of $P(f)=H^1_Q(M(f))$ is $Q(f)=Ext^3_S(M(f),S(-4))$, which is exactly the cokernel of $\delta$. It follows that
$$Q(f)=S/(a_1,a_2,a_3,a_4)(d+d_3-4).$$
In particular, $Q(f)$ is a finite dimensional vector space if and only if  $a_1,a_2,a_3,a_4$ is a regular sequence in $S$. 

 This completes the proof of the second  claim.
\endproof

\begin{ex}
\label{exNF4} 

(i) The surface $D:f=x^4-xyw^2+zw^3=0$ from Example \ref{exNF} is nearly free, and the generating syzygies are the following
$$r_1=(0,2y,3z,-w),  ~ r_2=(0,w,x,0),  ~ r_3=(w^2,4x^2,yw,0),  ~  r_4=(yw,6xz,y^2,2x^2).$$
It follows that $a_1=2x^2$, $a_2=-6xz$, $a_3=-y$ and $a_4=w$, and hence 
$$I_A=(x^2, xz,y,w).$$
In this case
 $a_1,a_2,a_3,a_4$ is not a regular sequence in $S$ and the proof above implies that in such a case
$$H^1_Q(M(f))_k \ne 0$$
for any integer $k \leq d+d_3-4$. As a matter of fact, all the examples of nearly free surfaces listed above have a similar property.

\noindent (ii) The examples with $H^1_Q(M(f))$ is a finite dimensional $\C$-vector space seem to be very rare. The following examples were kindly provided to us by Aldo Conca:
$$D_6:f=x^6 + x^4y^2 + y^5z + x^2y^3w=0 \text{ with exponents } (2,2,2),$$
$$ \  \  \  \   \   \  \  \  \  \  \ D_7:f= y^7 + xy^4w^2 + y^5w^2 + y^3w^4 + zw^6=0  \text{ with exponents } (2,2,3),$$
$$ \  \  \  \  D_8:f= x^8 + x^7z + x^3z^3w^2 + yz^4w^3=0 \text{ with exponents } (2,3,3),$$
and the two index family
$$D_{a,b}:f=x^{2a+2b-1}+x^{a+b-1} y^a z^b + y^{2a-1} z^{2b-1} w=0$$
with $a>1$, $b>1$ and
$a+b>4$.
In all these examples, it is easy to check that the surface is nearly free and the sequence 
$a_1,a_2,a_3,a_4$ is regular. 
For instance, for the surface $D_6$ one has
$(a_1,a_2,a_3,a_4)=(-15x, 15y,-15z+2w,9w),$ and hence
$$I_A=(x,y,z,w)=Q,$$
the maximal ideal in $S$.

\end{ex}

\begin{rk}
\label{locfree}
Consider the reflexive sheaf $Der(-log D)$ of logarithmic vector fields along $D$, which is the coherent sheaf on $\PP^3$ associated to the graded $S$-module $AR(f)(1)$.
By definition, a surface $D:f=0$ is free if and only if $Der(-log D)$ splits as a direct sum of line bundles. On the other hand,
Proposition \ref{prop1} (i) implies that if $H^0(\Sigma, \OO_{\Sigma}(k))=H^1_Q(M(f))_k \ne 0$ for $k<<0$, then  $Der(-log D)$ is not locally free by Lemma 3.2 in \cite{AY}.
Indeed, obvious exact sequences imply that
$$H^2(\PP^3, Der(-log D)(k))=H^1(\PP^3,\J_f(d+k))=H^0(\Sigma, \OO_{\Sigma}(d+k)),$$
where $\J_f$ is the ideal sheaf associated to $J_f$. In particular, for  all the examples of nearly free surfaces listed in this note except those in Example \ref{exNF4} (ii), the sheaf $Der(-log D)$ is not locally free.

On the other hand, again by Lemma 3.2 in \cite{AY}, for the nearly free surfaces listed in  Example \ref{exNF4} (ii), the sheaf $Der(-log D)$ is  locally free and not a direct sum of line bundles. Examples of sheaves $Der(-log D)$ which are locally free and not a direct sum of line bundles come also from the geometry, using the locally free arrangements, see \cite{MuSch, Yu} or \cite[Theorem 8.5]{DHA}. For instance, the arrangement in $\C^4$ given by
$$\A:f=\prod_{a=(a_0,a_1,a_2,a_3)}(a_0x+a_1y+a_2z+a_3w)=0$$
where $a=(a_0,a_1,a_2,a_3) \in \{0,1\}^4$, $a \ne (0,0,0,0)$, is locally free but not free,
see \cite[Example 4.5]{MuSch}. This arrangement,  going back to \cite{ER}, is nearly free with $d=15$ and $d_1=d_2=d_3=15$, as follows from the free resolution of $D^1_0=AR(f)$ given in \cite[Example 4.5]{MuSch}. Note that this arrangement is not tame in the sense of hyperplane arrangement theory, see for instance \cite{OT95, WY}, but it is tame as a surface in $\PP^4$ in the sense of Remark \ref{rktame}. The last claim follows by a direct computation using the Singular software \cite{Sing}.
\end{rk}

\begin{rk}
\label{rksymmetry}
It is interesting to note that the free and nearly surfaces in our examples, except the discriminant section $\D_4'$,  have all rather large symmetry groups. For instance, the surfaces  $D'_d$  introduced in Example \ref{exFREE2} (ii) 
admit all an effective $\C$-action given by 
$$t \cdot (x,y,z,w)=(x, y, z+ty, w-tx)$$
whose fixed point set is exactly the singular set $x=y=0$. This is related to the fact that in all these examples, except that in Proposition \ref{propD4'}, one has $d_1=1$. 

\end{rk}

\end{document}